\documentclass{elsarticle}

\usepackage{atbegshi}                          
\AtBeginDocument{\AtBeginShipoutNext{\AtBeginShipoutDiscard}}

\usepackage{amssymb}
\usepackage{amsmath}

\renewcommand{\thefootnote}{}

\newcommand{\NN}{{\mathbb N}}

\newtheorem{theorem}{Theorem}
\newtheorem{statement}{Statement}

\newtheorem{remark}[theorem]{Remark}

\newfont{\eurorm}{eurm10 scaled 1100}
\newfont{\eurorms}{eurm10 scaled 800}

\journal{arXiv.org}

\date{February 6, 2018}

\begin{document}

\begin{frontmatter}

\title{\bf \large A NEW METHOD FOR PROVING SOME INEQUALITIES RELATED TO SEVERAL SPECIAL FUNCTIONS}

\author{Tatjana Lutovac}
\ead{tatjana.lutovac@etf.bg.ac.rs}
\author{Branko Male\v sevi\' c${}^{\mbox{\scriptsize $\ast$}}$}
\ead{branko.malesevic@etf.bg.ac.rs}
\author{Marija Ra\v sajski}
\ead{marija.rasajski@etf.bg.ac.rs}

\address{\normalsize \rm School of Electrical Engineering, University of  Belgrade,   \\
Bulevar Kralja Aleksandra 73, 11000 Belgrade, Serbia}

\vspace*{-5.0 mm}

\begin{abstract}
In this paper we present a new approach to  proving some exponential inequalities involving the sinc function.  Power series expansions
are used to generate new polynomial inequalities that are sufficient to prove the given exponential inequalities.
\end{abstract}


\end{frontmatter}

\section{Introduction and preliminaries}

\renewcommand{\thefootnote}{}

\footnotetext{$\!\!\!\!\!\!{}^{\ast}\,$Corresponding author.                                                                   \\
Research of the first and second and  third author was supported in part by the Serbian Ministry of Education,
Science and Technological Development, under  Projects TR 32023,  ON 174032 \& III 44006 and ON 174033, respectively.}

\address{\rm School of Electrical Engineering, University of  Belgrade,   \\
Bulevar Kralja Aleksandra 73, 11000 Belgrade, Serbia}

Inequalities involving trigonometric and inverse trigonometric functions play an important role
and have many applications in science and engineering, see \cite{DSM_1970},
\cite{Rahmatollahi_DeAbreu_2012}$-$\cite{Cloud_Drachman_Lebedev_2014}, 
\cite{Malesevic_Rasajski_Lutovac_2017a}, \cite{Malesevic_Rasajski_Lutovac_2017b}. Of special interest are
inequalities with sinc function, i.e. $\displaystyle \mbox{\rm sinc}\,x\!=\!\mbox{\small $\dfrac{\sin x}{x}$}$
\mbox{${\big (}\displaystyle 0 \!<\! x \!\leq\! \mbox{\small $\dfrac{\pi}{2}$}{\big )}$.}
It is well-known that the sinc function is often used in signal processing, optics, radio transmission, sound recording, etc.\\

Starting from  {\sc Jordan}'s inequality \cite{DSM_1970},
\[
\displaystyle
\frac{2}{\pi}\leq\frac{\sin x}{x}\leq 1, ~~~~0<x \leq \frac{\pi}{2},
\]
and continuing with the polynomial bounds \cite{Qi_1996}, \cite{Zhu_2006}, \cite{Zhu_2006_II},
\cite{Agarwal_Kim_Sen_2009}, \cite{Debnath_Mortici_Zhu_2015}, some exponential bounds have recently been considered
\cite{Nishizawa_2017}.

In \cite{Nishizawa_2015} and \cite{Chen_Shi_Wang_Xiang_2017} the inequalities of the following form were studied:
\begin{equation}
\label{Ineq_1}
\left(1 - \mbox{\small $\displaystyle\frac{4(\pi-2)}{\pi^3}$}\,x^2\right)^{\!\alpha_1 x^3 + \beta_1 x^2 + \gamma_1 x + \delta_1}
>
\frac{\sin x}{x}
\end{equation}
and
\begin{equation}
\label{Ineq_2}
\left(1 - \mbox{\small $\displaystyle\frac{4(\pi-2)}{\pi^3}$}\,x^2\right)^{\!\alpha_2 x^3 + \beta_2 x^2 + \gamma_2 x + \delta_2}
<
\frac{\sin x}{x}
,
\end{equation}
for  $x \!\in\! \left(0,\mbox{\small $\displaystyle\frac{\pi}{2}$}\right)$ where  $\alpha_{j}, \beta_{j}, \gamma_{j}, \delta_{j}$ ($j\!=\!1,2$) are specified real coefficients.

\medskip
In this paper we present a new approach to proving  exponential inequalities of the above form. Using  the power series expansions of
the corresponding  functions and some newly developed approximation  techniques, we  reduce  exponential inequalities to the corresponding
polynomial inequalities that are  more  easily analysed and proved.

\smallskip
The application of our method is illustrated  in the proofs of the inequalities (\ref{Ineq_1}) and (\ref{Ineq_2}),
where the  coefficients  $\alpha_{j}, \beta_{j}, \gamma_{j}, \delta_{j}$ ($j\!=\!1,2$) are calculated    from the
constraints proposed  in \cite{Chen_Shi_Wang_Xiang_2017}.

\medskip
In the rest of this section we review some results that we use in our study.\\

Firstly, let us recall some well known  power series, such us:
\begin{equation}
\label{ln(1+x)}
\ln (1+x)
=
\displaystyle\sum\limits_{k=1}^{\infty}{(-1)^{k-1}\displaystyle\frac{x^k}{k}}, \qquad (-1 < x \leq 1)
\end{equation}
and
\begin{equation}
\label{ln(1-x)}
\ln (1-x)
=
-\!\displaystyle\sum\limits_{k=1}^{\infty}{\displaystyle\frac{x^k}{k}}, \qquad (-1 \leq x < 1).
\end{equation}

Also, in accordance with  \cite{Gradshteyn_Ryzhik_2015},  the following expansions hold:

\begin{equation}
\label{Series_ln_sin_x_over_x}
\ln \frac{\sin x}{x}
=
-\sum\limits_{k=1}^{\infty}{\frac{2^{2k-1}|B_{2k}|}{k(2k)!}x^{2k}}, \qquad (0 < x < \pi),
\end{equation}
\begin{equation}
\label{Series_ln_cos_x}
\ln \cos x
=
-\sum\limits_{k=1}^{\infty}{\frac{2^{2k-1}(2^{2k}-1)|B_{2k}|}{k(2k)!}x^{2k}}, \qquad (-\pi/2 < x < \pi/2),
\end{equation}
 where  $B_{i}$ ($i \!\in\! N$) are {\sc Bernoulli}'s numbers.

\medskip
The following Statement, which is a consequence of
Theorem 2 from   \cite{Wu_Debnath_2009},  was proved in \cite{Malesevic_Rasajski_Lutovac_2017c}.

\begin{statement}
\label{Statement_WD}
For the function $f\!:\!(a,b) \longrightarrow R$ let there exist the power series expansion$\,:$
\begin{equation}
f(x)
=
\displaystyle\sum_{k=0}^{\infty}{c_{k}(x-a)^k},
\end{equation}
for every $x \!\in\! (a,b)$, where $\{c_{k}\}_{k \!\in\! N_0}$ is the sequence of non-negative coefficients.
Then the following holds$\,:$
\begin{equation}
\label{WD_posledica}
\!
\begin{array}{c}
\displaystyle\sum_{k=0}^{n-1}{\!c_k(x-a)^k}
+
\frac{1}{(b-a)^n}
{\bigg (}
f(b-)
-
\displaystyle\sum_{k=0}^{n-1}{\!c_k(b-a)^k}
{\bigg )}(x-a)^{n}                                                                                            \\[1.5 ex]
\geq
f(x)
\geq
\displaystyle\sum_{k=0}^{n}{\!c_k(x-a)^{k}},
\end{array}
\end{equation}

\vspace*{-2.0 mm}

\noindent
for every $x \!\in\! (a,b)$ and $n \!\in\! N$.
\end{statement}
\begin{remark}
Note that in $(\ref{WD_posledica})$ strict inequalities hold except for some special cases of polynomial function.
\end{remark}

\section{Main results}

  We start this Section by deriving  some double-sided  inequalities
  that are consequences  of the power series expansions and  Statement 1,
  needed in the proofs of Theorems 1 and 2.

\medskip

\subsection{Some important double-sided inequalities}


\medskip
Using Statement \ref{Statement_WD} from (\ref{Series_ln_sin_x_over_x}) we get the following double-sided inequality:
\begin{equation}
\label{procena-ln-sin-x-over-x}
\begin{array}{l}
\;\;\;\;\,
-\displaystyle\sum\limits_{k=1}^{n}{
\mbox{\small $\displaystyle\frac{2^{2k-1}|B_{2k}|}{k(2k)!}$}x^{2k}}
\,>\,
\ln
\mbox{\small $\displaystyle\frac{\sin x}{x}$} \, >                                                                \\[3.0 ex]
\,>
-\displaystyle\sum\limits_{k=1}^{m-1}{
\mbox{\small $\displaystyle\frac{2^{2k-1}|B_{2k}|}{k(2k)!}$}x^{2k}}
+
{\Big (}\mbox{\small $\displaystyle\frac{2}{\pi}$}{\Big )}^{\!\!2m}\!\!
\left(
\ln \mbox{\small $\displaystyle\frac{2}{\pi}$}
-
\displaystyle\sum\limits_{k=1}^{m-1}{
\mbox{\small $\displaystyle\frac{2^{2k-1}|B_{2k}|}{k(2k)!}$}
{\Big (}\frac{\pi}{2}{\Big )}^{2k}}\right)
\!x^{2m}\!,
\end{array}
\end{equation}
for $x \!\in\! \left(0,\mbox{$\displaystyle\frac{\pi}{2}$}\right]$, where $n, m \!\in\! N$.

\medskip
Based on Statement \ref{Statement_WD} from (\ref{Series_ln_cos_x}) we have:
\begin{equation}
\label{Procena-ln-cos-x}
\!\!
\begin{array}{l}
\;\;\;
-\!\!\displaystyle\sum\limits_{k=1}^{n}{\!
\mbox{\small $\displaystyle\frac{2^{2k-1}(2^{2k}\!-\!1)|B_{2k}|}{k(2k)!}$}x^{2k}}
>
\ln
\cos x \, >                                                                                                   \\[3.0 ex]
\!>\!
-\!\!\displaystyle\sum\limits_{k=1}^{m-1}{\!
\mbox{\small $\displaystyle\frac{2^{2k-1}(2^{2k}\!-\!1)|B_{2k}|}{k(2k)!}$}x^{2k}}
\!+\!
{\Big (}\mbox{\small $\displaystyle\frac{1}{c}$}{\Big )}^{\!2m}\!\!
\left(\!
\ln \cos c
-
\!\!\displaystyle\sum\limits_{k=1}^{m-1}{
\mbox{\small $\displaystyle\frac{2^{2k-1}(2^{2k}\!-\!1)|B_{2k}|}{k(2k)!}$}
c^{2k}}\!\right)
\!x^{2m}\!,
\end{array}
\!\!\!\!\!\!\!\!\!\!\!\!\!\!
\!\!\!\!\!\!\!\!\!\!\!\!\!\!
\end{equation}
for $x \!\in\! \left(0,c\right]$, where $0 \!<\! c \!<\! \mbox{$\displaystyle\frac{\pi}{2}$}$, $n, m \!\in\! N$.

\medskip
Using {\sc Leibniz}'s theorem applied to (\ref{ln(1+x)}) we obtain the following double-sided inequality:
\begin{equation}
\label{Procena-ln(1+x)}
\displaystyle\sum\limits_{k=1}^{2\ell_1-1}{(-1)^{k-1}\displaystyle\frac{x^k}{k}}
\,>\,
\ln (1+x)
\,>\,
\displaystyle\sum\limits_{k=1}^{2\ell_2}{(-1)^{k-1}\displaystyle\frac{x^k}{k}},
\end{equation}
for $x \!\in\! (0,1)$ and  $\ell_1, \ell_2 \!\in\! N$.

\break

Based on Statement \ref{Statement_WD} from (\ref{ln(1-x)}) we get:
\begin{equation}
\label{procena-ln(1-x)}
-\displaystyle\sum\limits_{k=1}^{n}{\displaystyle\frac{x^{k}}{k}}
\,>\,
\ln (1-x)
\,>\,
-\displaystyle\sum\limits_{k=1}^{m-1}{\displaystyle\frac{x^{k}}{k}}
+
\frac{1}{c^{m}}\!\!\left(\!\ln(1-c) + \displaystyle\sum\limits_{k=1}^{m-1}{\displaystyle\frac{c^{k}}{k}}\!\right)\!x^{m},
\end{equation}
for $x \!\in\! (0,c)$, where $ 0<c<1$, and $n, m \!\in\! N$.

\medskip

\subsection{On some constraints and their consequences}
\label{About-constraints}

\medskip
Now we show that the constraints proposed in \cite{Chen_Shi_Wang_Xiang_2017}  induce linear relations
among the coefficients  $\alpha_{j}, \beta_{j}, \gamma_{j}, \delta_{j}$ of the polynomials in the exponents of the functions on the left-hand side in (\ref{Ineq_1}) and (\ref{Ineq_2}).

\medskip

Let us consider  function  $f$ defined by:
\begin{equation}
\label{Def_f(x)}
f(x)
\!=\!
\left(1 - \mbox{\small $\displaystyle\frac{4(\pi-2)}{\pi^3}$}\,x^2\right)^{\!\alpha x^3 + \beta x^2 + \gamma x + \delta}
\!-\;
\frac{\sin x}{x},
\end{equation}
for some real coefficients $\alpha, \beta, \gamma, \delta$,  where $x \!\in\! \left(0,\mbox{\small $\displaystyle\frac{\pi}{2}$}\right]$.\\

\medskip
Let us notice that for  function $f$ the following two conditions are always satisfied:

\begin{equation}
f(0+)=0
\qquad\mbox{i}\qquad
f'(0+)=0.
\end{equation}
Next, the following holds for  function $f$:

\begin{equation}
\label{fj''(0)=0}
f''(0+)=0
\;\Longleftrightarrow\;
\delta
=
\mbox{\small $\displaystyle\frac{\pi^3}{24(\pi-2)}$}.
\end{equation}
Also, for function $f$ we have:

\begin{equation}
\label{fj(Pi/2)=0}
f\!\left(\mbox{\small $\displaystyle\frac{\pi}{2}$}\right)=0
\;\Longleftrightarrow\;
\delta
\!=\!
-
\mbox{\small $\displaystyle\frac{\pi^3}{8}$}\alpha
-
\mbox{\small $\displaystyle\frac{\pi^2}{4}$}\beta
-
\mbox{\small $\displaystyle\frac{\pi}{2}$}\gamma.
\end{equation}
Adding the condition $f'\!\left(\mbox{\small $\displaystyle\frac{\pi}{2}$}\right)\!=\!0$ we get:
\begin{equation}
\label{fj(Pi/2)=0 & fj'(Pi/2)=0}
f\!\left(\mbox{\small $\displaystyle\frac{\pi}{2}$}\right)\!=\!0
\wedge
f'\!\left(\mbox{\small $\displaystyle\frac{\pi}{2}$}\right)\!=\!0
\Longleftrightarrow
\left\{\!\!\!
\begin{array}{c}
\delta
\!=\!
-
\mbox{\small $\displaystyle\frac{\pi^3}{8}$}\alpha
-
\mbox{\small $\displaystyle\frac{\pi^2}{4}$}\beta
-
\mbox{\small $\displaystyle\frac{\pi}{2}$}\gamma                                                              \\[2.5 ex]
\gamma
\!=\!
-
\mbox{\small
$\displaystyle\frac{\pi^2 (\ln \pi \!-\! \ln 2 )(3\,\pi\,\alpha \!+\! 4\,\beta) + 8 (\pi \!-\! 3)}{4\,\pi\,(\ln \pi \!-\! \ln 2 )}$}
\end{array}
\!\!\!\right\}\!.\!\!\!
\end{equation}
Further, adding the condition $f''\!\left(\mbox{\small $\displaystyle\frac{\pi}{2}$}\right)\!=\!0$ we get:
\begin{equation}
\label{fj(Pi/2)=0 & fj'(Pi/2)=0 & fj''(Pi/2)=0}
\begin{array}{rl}
&
f\!\left(\mbox{\small $\displaystyle\frac{\pi}{2}$}\right)\!=\!0
\wedge
f'\!\left(\mbox{\small $\displaystyle\frac{\pi}{2}$}\right)\!=\!0
\wedge
f''\!\left(\mbox{\small $\displaystyle\frac{\pi}{2}$}\right)\!=\!0                                        \\[2.0 ex]
\!\!\Longleftrightarrow\!\!
&
\left\{\!\!\!
\begin{array}{c}
\delta
\!=\!
\!-\!
\mbox{\small $\displaystyle\frac{\pi^3}{8}$}\alpha
\!-\!
\mbox{\small $\displaystyle\frac{\pi^2}{4}$}\beta
\!-\!
\mbox{\small $\displaystyle\frac{\pi}{2}$}\gamma                                                             \\[2.5 ex]
\gamma
\!=\!
\!-\!
\mbox{\small
$\displaystyle\frac{\pi^2 (\ln \pi \!-\! \ln 2 )(3\,\pi\,\alpha \!+\! 4\,\beta) + 8 (\pi \!-\! 3)}{4\,\pi\,(\ln \pi \!-\! \ln 2 )}$}
                                                                                                              \\[2.5 ex]
\beta
\!=\!
\mbox{\small
$\displaystyle\frac{-3\,\pi^3(\ln 2 \!-\! \ln \pi )^2 \alpha
+
(\pi\!-\!2)(3\pi\!-\!6)(\ln 2\!-\!\ln \pi)
+
8(\pi\!-\!3)}{
2 \pi^2 (\ln 2 \!-\! \ln \pi)^2}$}
\end{array}
\!\!\!\right\}\!.\!\!\!
\end{array}
\end{equation}

Let us note that the conditions connecting the coefficients $\alpha, \beta, \gamma, \delta$ are linear.

\medskip

\subsection{Proofs of exponential inequalities $(1)$ and $(2)$}
\label{Proofs}

\medskip

In this section we determine the real coefficients $\alpha_j, \beta_j, \gamma_j, \delta_j$, ($j=1,2$) from the constraints proposed
in \cite{Chen_Shi_Wang_Xiang_2017}, and in Theorem 1 and Theorem 2 we prove the corresponding inequalities.

\setcounter{theorem}{0}

\begin{theorem}
Let  function $f_1$ be defined in the interval {\small $\left(0,\mbox{\small $\displaystyle\frac{\pi}{2}$}\right]$} by

\begin{equation}
f_1(x)
=
\left(1 - \mbox{\small $\displaystyle\frac{4(\pi-2)}{\pi^3}$}\,x^2\right)^{\!\alpha_1 x^3 + \beta_1 x^2 + \gamma_1 x + \delta_1}
\!-\;
\frac{\sin x}{x}
\end{equation}
and let the following conditions hold$:$
\begin{equation}
\label{Uslovi_1}
f_{1}(0+)
\!=\!
f_{1}'(0+)
\!=\!
f_{1}''(0+)
\!=\!
0, \;
f_{1}\!\left(\mbox{\small $\displaystyle\frac{\pi}{2}$}\right)
\!=\!
f_{1}'\!\left(\mbox{\small $\displaystyle\frac{\pi}{2}$}\right)
\!=\!
f_{1}''\!\left(\mbox{\small $\displaystyle\frac{\pi}{2}$}\right)
\!=\!
0.
\end{equation}
Then:
\begin{equation}
\label{Koeficijenti_1}
\begin{array}{rcl}
\alpha_1
\!\!\!&\!\!=\!\!&\!\!\!
\mbox{\footnotesize $\displaystyle\frac{
\left(-\pi^{3}\!+\!24\pi\!-\!48 \right)\!\ln^2\! \frac{\pi}{2}
\!-\!
3\left(\pi\!-\!2\right)\!\left(3\pi^{2}\!-\!20\pi\!+\!36\right)
\ln \! \frac{\pi}{2}
+
24\!\left(\pi\!-\!3\right)\!\left(\pi\!-\!2\right)^{2}
}{
3
\left(\pi - 2 \right)\pi^{3}\ln^{2}\! \frac{\pi}{2}
}$},                                                                                                           \\[1.5 ex]
\beta_1
\!\!\!&\!\!=\!\!&\!\!\!
\mbox{\footnotesize $\displaystyle\frac{(\pi^{3}-24\pi+48)\ln^2 \!\frac{\pi}{2}
+
6(\pi-2)(\pi-4)^{2} \ln \!\frac{\pi}{2}
-
16(\pi-3)(\pi-2)^{2}
}{
2(\pi-2)\pi^{2}\ln^2\! \frac{\pi}{2}}$},                                                                      \\[2.0 ex]
\gamma_1
\!\!\!&\!\!=\!\!&\!\!\!
\mbox{\footnotesize $\displaystyle\frac{(-\pi^{3}\!+\!24\pi\!-\!48)\ln^2\! \frac{\pi}{2}
+
(3\pi\!-\!10)(\pi\!-\!6)(\pi\!-\!2) \ln \!\frac{\pi}{2}
+
8(\pi\!-\!3)(\pi\!-\!2)^{2}
}{
4\pi(\pi\!-\!2)\ln^{2} \!\frac{\pi}{2}
}$},                                                                                                           \\[1.0 ex]
\delta_1
\!\!\!&\!\!=\!\!&\!\!\!
\displaystyle\frac{\pi^3}{24(\pi-2)}
\end{array}
\end{equation}
and
\begin{equation}
\label{Nejednakost_1}
f_1(x) > 0,
\end{equation}
for every $x \!\in\! \left(0,\mbox{\small $\displaystyle\frac{\pi}{2}$}\right]$.
\end{theorem}
{\bf Proof.}
As discussed in Subsection~\ref{About-constraints}, the system derived from the conditions in (\ref{Uslovi_1}) can be reduced to the system of linear algebraic equations from  (\ref{fj''(0)=0}) and (\ref{fj(Pi/2)=0 & fj'(Pi/2)=0 & fj''(Pi/2)=0}), in variables $\alpha_1$, $\beta_1$, $\gamma_1$ i $\delta_1$.
The solution to this system, i.e. the coefficients of the polynomial in the exponent of the function $f_1$, are shown in (\ref{Koeficijenti_1}). The corresponding numerical values are: $\alpha_1 =  -0.0277933961 \ldots, \beta_1 = 0.0136111520\ldots, \gamma_1 = -0.0366389131\ldots$, $\delta_1 = 1.13168930\ldots\;.$

\medskip
\noindent


We consider  inequality (\ref{Nejednakost_1}) in its equivalent form:
\begin{equation}
\label{EKV_Nejednakost_1}
\left( \alpha_1 x^3 + \beta_1 x^2 + \gamma_1 x + \delta_1 \right)
\ln\!\left(1 - \mbox{\small $\displaystyle\frac{4(\pi-2)}{\pi^3}$}\,x^2\right)
>
\ln\mbox{\small $\displaystyle\frac{\sin x}{x}$},
\end{equation}
for $x \!\in\! \left(0,\mbox{\small $\displaystyle\dfrac{\pi}{2}$}\right]$. \\

We will prove that the corresponding function of the difference$:$

\begin{equation}
g(x)
=
\left( \alpha_1 x^3 + \beta_1 x^2 + \gamma_1 x + \delta_1 \right)
\ln\!\left(1 - \mbox{\small $\displaystyle\frac{4(\pi-2)}{\pi^3}$}\,x^2\right)
-
\ln \mbox{\small $\displaystyle\frac{\sin x}{x}$}
\end{equation}
is positive for every $x \!\in\! \left(0,\mbox{\small $\displaystyle\frac{\pi}{2}$}\right]$.

\smallskip
It is easy to verify that the polynomial of the third degree
$$
P_{3}(x) = \alpha_1 x^3 + \beta_1 x^2 + \gamma_1 x + \delta_1
$$
is positive in the interval $\left(0,\mbox{\small $\displaystyle\frac{\pi}{2}$}\right]$.\\

The proof of positivity of  function $g(x)$ is divided into two parts. In the first part we will consider the function around zero, and in the second part we focus on the remaining part of the interval  $\left(0,\mbox{\small $\displaystyle\frac{\pi}{2}$}\right]$.

{\bf 1.} Let us notice that
\begin{equation}
\mbox{\small $\displaystyle\frac{4(\pi-2)}{\pi^3}$}\,x^2 \leq \frac{\pi-2}{\pi} < 1,
\end{equation}
for $x \!\in\! \left(0,\mbox{\small $\displaystyle\frac{\pi}{2}$}\right]$. Therefore, we can use  inequality (\ref{procena-ln(1-x)}) in the form:
\begin{equation}
\begin{array}{l}
\ln\!\left(\!1 \!-\! \mbox{\small $\displaystyle\frac{4(\pi-2)}{\pi^3}$}x^2\!\right) >                              \\[2.0 ex]
>
-\!
\displaystyle\sum_{k=1}^{n-1}{\mbox{\small $\displaystyle\frac{4^{k}(\pi\!-\!2)^{k}}{\pi^{3k}k}$}x^{2k}}
+
\left(\mbox{\small $\displaystyle\frac{2}{\pi}$}\right)^{\!2n}
\!\!
\left(\!
\ln\! \mbox{\small $\displaystyle\frac{2}{\pi}$}
+\!
\displaystyle\sum_{k=1}^{n-1}{\mbox{\small $\displaystyle\frac{4^{k}(\pi\!-\!2)^{k}}{\pi^{3k}k}$}
\left(\!\mbox{\small $\displaystyle\frac{\pi}{2}$}\!\right)^{2k}}
\!\right)\!x^{2n},
\end{array}
\end{equation}
for $x \!\in\! \left(0,\mbox{\small $\displaystyle\frac{\pi}{2}$}\right]$ i $n \!\in\! N$, $n \geq 2$.
Also, from (\ref{procena-ln-sin-x-over-x}) the following holds$:$
\begin{equation}
-
\ln\!\mbox{\small $\displaystyle\frac{\sin x}{x}$}
>
\displaystyle\sum\limits_{k=1}^{m}{
\mbox{\small $\displaystyle\frac{2^{2k-1}|B_{2k}|}{k(2k)!}$}x^{2k}},
\end{equation}
for $x \!\in\! \left(0,\mbox{\small $\displaystyle\frac{\pi}{2}$}\right]$ and $m \!\in\! N$.

\smallskip

Now let us construct the approximating polynomial:
\begin{equation}
\!\!\!\!
\begin{array}{rcl}
G_{n,m}(x)
\!\!\!&\!\!=\!\!&\!\!\!
P_3(x) \!\left(
\!-\!
\displaystyle\sum_{k=1}^{n-1}{\mbox{\small $\displaystyle\frac{4^{k}(\pi\!-\!2)^{k}}{\pi^{3k}k}$}\,x^{2k}}
\!+\!
\left(
\mbox{\small $\displaystyle\frac{2}{\pi}$}\right)^{\!\!2n}
\!\!
\left(\!
\ln\! \mbox{\small $\displaystyle\frac{2}{\pi}$}\!
+\!\!
\displaystyle\sum_{k=1}^{n-1}{\mbox{\small $\displaystyle\frac{4^{k}(\pi\!-\!2)^{k}}{\pi^{3k}k}$}
\left(\!\mbox{\small $\displaystyle\frac{\pi}{2}$}\!\right)^{\!2k}}
\!\right)\!x^{2n}\!\!\right)                                                                                  \\[3.5 ex]
\!\!\!&\!\! \!\!&\!\!\!
+\,
\displaystyle\sum\limits_{k=1}^{m}{
\mbox{\small $\displaystyle\frac{2^{2k-1}|B_{2k}|}{k(2k)!}$}x^{2k}},
\end{array}\!\!\!\!\!\!\!\!\!\!\!\!\!\!\!\!\!\!\!\!\!\!\!\!\!\!\!\!\!\!
\end{equation}
for $x \!\in\! \left(0,\mbox{\small $\displaystyle\frac{\pi}{2}$}\right]$. \\
 We conclude that the following inequality  holds:
\begin{equation}
\label{G-polynomial-C1}
\displaystyle
g(x)>G_{n,m}(x),
\end{equation}
for $x \!\in\! \left(0,\mbox{\small $\displaystyle\frac{\pi}{2}$}\right]$ and  $m,n \!\in\! N,  \, \, n\geq 2$.

If we specify the values $m,n \!\in\! N$, we can determine the zeros and the sign of the polynomial $G_{n,m}$. For example, for $n\!=\!3$ and $m\!=\!2$ it is easy to verify that $G_{3,2}(x)>0$ for every $x \!\in\! (0, c_1)$, where $c_1$ is a positive number smaller than the first positive root of the considered polynomial, and that is $x =0.925930 \ldots\;$. Let us set: $c_1=0.92$. Let us notice that the polynomial $G_{3,2}$ is of degree $9$.

As   $\displaystyle g(x)>G_{3,2}(x)$ holds
for every $x \!\in\! \left(0,\mbox{\small $\displaystyle\frac{\pi}{2}$}\right]$, we conclude that $g(x)>0$ for every $x \!\in\! (0, c_1)$, i.e.  inequality (\ref{Nejednakost_1}) holds in the interval $(0, c_1)$.\\

\medskip
{\bf 2.}  To prove that inequality (\ref{Nejednakost_1}) holds in the remaining part of the interval $\left(0,\mbox{\small $\displaystyle\frac{\pi}{2}$}\right]$, we will introduce a change of variable $x\!=\!\mbox{\small $\displaystyle\frac{\pi}{2}$}-t$ in  inequality (\ref{EKV_Nejednakost_1}), $t \!\in\! \left[0,\mbox{\small $\displaystyle\frac{\pi}{2}$}\right)$. It is sufficient to prove that the inequality holds for \mbox{ $t \!\in\! \left[0,c_2\right]$}, where { $c_2 \!=\! 0.65 \!>\! \mbox{\small $\displaystyle\frac{\pi}{2}$} \!-\! 0.92
\!=\! 0.650796327\ldots\,$}. \\

Now let us consider the equivalent inequality:
\begin{equation}
\label{Conj._1_Ineq._3}
Q_{3}(t)
\ln\!\left(\!
1
-
\mbox{\small $\displaystyle\frac{4(\pi\!-\!2)}{\pi^3}$}\,\!\left(\!\mbox{\small $\displaystyle\frac{\pi}{2}$}\!-\!t\!\right)^{\!2}
\right)
>
\ln\! \mbox{\small $\displaystyle\frac{\sin \!\left(\!\mbox{\small $\displaystyle\frac{\pi}{2}$}\!-\!t\!\right)\!}{
\!\!\mbox{\small $\displaystyle\frac{\pi}{2}$}\!-\!t\!\!}$} \!
=
\ln \cos t
-
\ln\!\left(\mbox{\small $\displaystyle\frac{\pi}{2}$}\!-\!t\right)\!,
\end{equation}
for $t \!\in\! \left[0,c_2\right]$, where $Q_{3}(t) = P_{3}\!\left(\mbox{\small $\displaystyle\frac{\pi}{2}$}-t\right)$. \\

We will prove that the corresponding function of the difference:
\begin{equation}
h(t)
=
Q_{3}(t)
\ln\!\left(\!
1
\!-\!
\mbox{\small $\displaystyle\frac{4(\pi\!-\!2)}{\pi^3}$}\,\!\left(\mbox{\small $\displaystyle\frac{\pi}{2}$}\!-\!t\right)^{\!2}
\right)
-
\ln \cos t
+
\ln\!\left(\mbox{\small $\displaystyle\frac{\pi}{2}$}-t\right)\!,
\end{equation}
is positive for every $t \!\in\! \left[0, c_2\right]$.\\

The positivity of the polynomial $P_3(x)$ in the interval $\left(0,\mbox{\small $\displaystyle\frac{\pi}{2}$}\right]$ yields positivity of the polynomial
\mbox{$Q_{3}(t)$}.

Since the following holds:

\begin{equation}
\ln\!\left(1-\mbox{\small $\displaystyle\frac{4(\pi-2)}{\pi^3}$}\,\!\left(\!\mbox{\small $\displaystyle\frac{\pi}{2}$}-t\!\right)^{\!2}\right)
=
\ln \! \mbox{\small $\displaystyle\frac{2}{\pi}$}
+
\ln \! \left(\! 1 + \mbox{\small $\displaystyle\frac{2(\pi-2)t(\pi-t)}{\pi^2}$}\!\right)
\end{equation}
and
\begin{equation}
\frac{2(\pi-2)t(\pi-t)}{\pi^2} < 1,
\end{equation}
for every $t \!\in\! \left(0,c_2\right]$, we can use  inequality (\ref{procena-ln(1-x)}) in the following form:
\begin{equation}
\ln\!\left(1-\mbox{\small $\displaystyle\frac{4(\pi-2)}{\pi^3}$}\,\!\left(\!\mbox{\small $\displaystyle\frac{\pi}{2}$}-t\!\right)^{\!2}\right)
\geq \,
\ln \! \mbox{\small $\displaystyle\frac{2}{\pi}$}
+
\displaystyle\sum_{k=1}^{n}{(-1)^{k-1}\mbox{\small $\displaystyle\frac{2^k(\pi-2)^k}{k \pi^{2k}}$}\,t^k(\pi-t)^k}\!,
\end{equation}
for  $t \!\in\! \left[0,c_2\right]$ i $n=2 \ell \!\in\! N$.\\

Next, we can use  inequality (\ref{Procena-ln-cos-x})  in the following form:
\begin{equation}
-
\ln \cos t
\geq
\displaystyle\sum\limits_{k=1}^{m_1}{
\,\mbox{\small $\displaystyle\frac{2^{2k-1}(2^{2k}\!-\!1)|B_{2k}|}{k(2k)!}$}\,t^{2k}}
\end{equation}
for $t \!\in\! (0,c_2]$ and $m_1 \!\in\! N$. \\

Finally, let us notice that the following holds true:
\begin{equation}
\ln\!\left(\mbox{\small $\displaystyle\frac{\pi}{2}$}-t\right)
=
\ln\!\mbox{\small $\displaystyle\frac{\pi}{2}$}
+
\ln\!\left(1-\mbox{\small $\displaystyle\frac{2t}{\pi}$}\right)
\end{equation}
and
\begin{equation}
\frac{2t}{\pi} < 1,
\end{equation}
for $t \!\in\! \left[0,c_2\right]$.

\medskip
 We conclude that:
\begin{equation}
\!\!\!\!
\begin{array}{rcl}
\ln\!\left(\mbox{\small $\displaystyle\frac{\pi}{2}$}-t\right)
\!\!&\!\!=\!\!&\!\!
\ln \mbox{\small $\displaystyle\frac{\pi}{2}$}
-
\ln\!\left(1-\mbox{\small $\displaystyle\frac{2t}{\pi}$}\right)\geq                                               \\[2.0 ex]
\!\!&\!\!\geq \!\!&\!\!
\ln \mbox{\small $\displaystyle\frac{\pi}{2}$}
-\!
\displaystyle\sum_{k=1}^{m_2-1}{\!\!\mbox{\small $\displaystyle\frac{2^k}{\pi^kk}$}\,t^k\!}
+
\left(\mbox{\small $\displaystyle\frac{1}{c_2}$}\right)^{\!m_2}
\!\!
\left(
\!
\ln\!\left(1\!-\!\mbox{\small $\displaystyle\frac{2c_2}{\pi}$}\right)
\!+\!\!
\displaystyle\sum_{k=1}^{m_2-1}{\!\mbox{\small $\displaystyle\frac{2^kc_2^k}{\pi^kk}$}}
\right)\!t^{m_2}\,,
\end{array}
\!\!\!
\end{equation}
for $t \!\in\! \left[0,c_2\right]$ and $m_2 \!\in\! N$.

\smallskip
\smallskip
Let us now consider the approximating polynomial:
\begin{equation}
\!\!\!\!
\begin{array}{rcl}
H_{n,m_1,m_2}(t)
\!\!&\!\!=\!\!&\!\!
Q_{3}(t) \!\cdot\!\! \left(
\ln \mbox{\small $\displaystyle\frac{2}{\pi}$}
+
\displaystyle\sum_{k=1}^{n}{(-1)^{k-1}\mbox{\small $\displaystyle\frac{2^k(\pi-2)^k}{k \pi^{2k}}$}\,t^k(\pi-t)^k}
\right)                                                                                                       \\[1.5 ex]
\!\!&\!\! \!\!&\!\!
+
\displaystyle\sum\limits_{k=1}^{m_1}{
\mbox{\small $\displaystyle\frac{2^{2k-1}(2^{2k}\!-\!1)|B_{2k}|}{k(2k)!}$}t^{2k}}                             \\[1.5 ex]
\!\!&\!\! \!\!&\!\!
+
\ln \mbox{\small $\displaystyle\frac{\pi}{2}$}
-\!\!
\displaystyle\sum_{k=1}^{m_2-1}{\!\!\mbox{\small $\displaystyle\frac{2^k}{\pi^kk}$}\,t^k\!}
+
\left(\mbox{\small $\displaystyle\frac{1}{c_2}$}\right)^{\!m_2}
\!\!
\left(\!
\ln\!\left(1\!-\!\mbox{\small $\displaystyle\frac{2c_2}{\pi}$}\right)
\!+\!\!
\displaystyle\sum_{k=1}^{m_2-1}{\!\mbox{\small $\displaystyle\frac{2^kc_2^k}{\pi^kk}$}}
\!\right)\!t^{m_2},
\end{array}
\!\!\!\!\!\!
\end{equation}
for $t \!\in\! \left[0, c_2\right]$, $m,n \!\in\! N$.

\break

As before, for the specified values  of \mbox{$n\!=\!6$}, $m_1\!=\!4$ and $m_2\!=\!8$, we can determine the zeros
and the sign of the above-mentioned polynomial (for example applying  {\sc Sturm}'s theorem). Let us notice that
the degree of the polynomial $H_{6,4,8}(t)$ is equal to $15$ and that $t = 0.789165 \ldots\;$ is the first positive
root of the polynomial $H_{6,4,8}(x)$. \\

Therefore, \[h(t) > H_{6,4,8}(t) > 0,\] for every $t \!\in\! \left[0, c_2\right]$.

\smallskip
This proves  inequality (\ref{Conj._1_Ineq._3}) in the interval $\left[0,c_2\right]$, and  hence
$g(x)>0$ holds for every $x \!\in\! \left[\mbox{\small $\displaystyle\frac{\pi}{2}$}-c_2,\mbox{\small $\displaystyle\frac{\pi}{2}$}\right]$.\\

Finally, it follows from {\bf 1.}$\;$and {\bf 2.}  that inequality  (\ref{Nejednakost_1})
holds for every \mbox{$x \!\in\! \left(0,\mbox{\small $\displaystyle\frac{\pi}{2}$}\right].$}
\hfill
$\Box$

\title{\bf \large Conjecture 2 }

\medskip
\medskip
\noindent
\begin{theorem}
Let the function
$$
\displaystyle
 f_2(x) \! \!= \!\!
\left(\!1 - \mbox{\small $\displaystyle\frac{4(\pi-2)}{\pi^3}$}\,x^2\!\right)^{\alpha_2 x^3 + \beta_2 x^2  +
\delta_2}
\! -  \frac{\sin x}{x},
$$
for  $x \!\in\! \left(0,\mbox{\small $\displaystyle\frac{\pi}{2}$}\right]$ satisfy the following conditions$:$
\begin{equation}
\label{Uslovi_2}
f_{2}(0+)
\!=\!
f_{2}'(0+)
\!=\!
f_{2}''(0+)
\!=\!
0,\,\,\,\,
f_{2}\!\left(\mbox{\small $\displaystyle\frac{\pi}{2}$}\right)
\!=\!
f_{2}'\!\left(\mbox{\small $\displaystyle\frac{\pi}{2}$}\right)
\!=\!
0.
\end{equation}
Then:
\begin{equation}
\label{Koeficijenti_2}
\begin{array}{rcl}\displaystyle
\alpha_2
\!\!&\!\!=\!\!&\!\!
\mbox{\footnotesize $\displaystyle
- \frac{2}{3}
\frac{12(\pi-2)(\pi-3)\, - \,
(48-24\pi+ \pi^{3})\ln \frac{\pi}{2}
}
{
\pi^{3} (\pi-2) \ln\frac{\pi}{2} \,
}
$},                                                                                           \\[2.2 ex]
\beta_2
\!\!&\!\!=\!\!&\!\!
\mbox{\small $\displaystyle
\frac{8(\pi-2)(\pi-3) \, - \,(48-24\pi+ \pi^{3})\ln \frac{\pi}{2}\,}
     {2 \pi^{2} (\pi-2) \ln \frac{\pi}{2}}$
},                                                                                        \\[2.2 ex]
\delta_2
\!\!&\!\!=\!\!&\!\!
\mbox{\small $\displaystyle\frac{\pi^3}{24(\pi-2)}$}
\\[2.0 ex]
\end{array}
\end{equation}
and
\begin{equation}\label{C2-ineq}
f_2(x)
<
0,
\end{equation}
for every $x \!\in\! \left(0,\mbox{\small $\displaystyle\frac{\pi}{2}$}\right]$.
\end{theorem}


~\\
\noindent {\bf Proof.}
As discussed in  Subsection~\ref{About-constraints},
the conditions in (\ref{Uslovi_2}) yield a  system   of linear equations  {\big (}shown in (\ref{fj''(0)=0}) and
(\ref{fj(Pi/2)=0 & fj'(Pi/2)=0}){\big )}  in variables $\, \alpha_2, \beta_2$ and  $ \delta_2$.
 The symbolic values (\ref{Koeficijenti_2}) of  $\alpha_2, \beta_2$ and $ \delta_2$
are obtained by solving this system.

Notice that it is easy to get  numeric values:
$\alpha_2 =  -0.0129442047 \ldots$ ,
$\beta_2 = -0.0330389552 \ldots$,
and $\delta_2 =  1.13168930\ldots.$

The exponential inequality~(\ref{C2-ineq}) is equivalent to the following inequality:
\begin{equation}
\label{ekviv-C2-ineq}
F_2(x) < 0
\end{equation}
for every $x \!\in\! \left(0,\mbox{\small $\displaystyle\frac{\pi}{2}$}\right]$, where
$$
F_2(x)= \left(\alpha_2 x^3 + \beta_2 x^2 +  \delta_2\right)
\ln\!\left(\!1 - \mbox{\small $\displaystyle\frac{4(\pi-2)}{\pi^3}$}\,x^2\!\right)
-
\ln\! \mbox{\small $\displaystyle\frac{\sin x}{x}$}.
$$
It is not difficult  to check that
\begin{equation}
\label{Polinom-P3}
P_{3}(x) = \alpha_2 x^3 + \beta_2 x^2  +
\delta_2 >0
\end{equation}
for every $x \!\in\! \left(0,\mbox{\small $\displaystyle\frac{\pi}{2}$}\right]$.

~\\
Based on   (\ref{procena-ln-sin-x-over-x}) and (\ref{procena-ln(1-x)}), for $x \!\in\! \left(0,\mbox{\small $\displaystyle\frac{\pi}{2}$}\right]$ and  $\, n,m \!\in\! N, \, n\geq 1, \, m\geq 2 $
the following inequalities hold:
$$
\begin{array}{l}
-\ln\!
\mbox{\small $\displaystyle\frac{\sin x}{x}$}
 \, < \,
\displaystyle\sum\limits_{k=1}^{m-1}{
\mbox{\small $\displaystyle\frac{2^{2k-1}|B_{2k}|}{k(2k)!}$}x^{2k}}
-{\Big (}\mbox{\small $\displaystyle \frac{2}{\pi}$}{\Big )}^{\!\!2m}\!\!
\left(\ln \frac{2}{\pi} \, - \,
\displaystyle\sum\limits_{k=1}^{m-1}{
\mbox{\small $\displaystyle\frac{2^{2k-1}|B_{2k}|}{k(2k)!}$}
\left(\frac{\pi}{2}\right)^{2k}}\right)
\!x^{2m}
\end{array}
$$
and
$$
\begin{array}{l}
\ln\!\left(1 - \mbox{\small $\displaystyle\frac{4(\pi-2)}{\pi^3}$}\,x^2\right) \;
< \; -
\displaystyle\sum_{k=1}^{n}{\mbox{\small $\displaystyle\frac{4^{k}(\pi\!-\!2)^{k}}{\pi^{3k}k}$}\,x^{2k}}.
\end{array}
$$

\smallskip
Finally, for $x\in \left(0,\mbox{\small $\displaystyle\frac{\pi}{2}$}\right]$, we have:
$$
\!\!\!\!\!\!\!\!
\begin{array}{lc}
F_{2}(x)
\!\! \!\!\! &
< \! - P_{3}(x)\!
\displaystyle\sum_{k=1}^{n}\!{\mbox{\small $\displaystyle\frac{4^{k}(\pi\!-\!2)^{k}}{\pi^{3k}k}$}x^{2k}}
\!+\!
\displaystyle\sum\limits_{k=1}^{m-1}\!{\mbox{\small $\displaystyle\frac{2^{2k-1}|B_{2k}|}{k(2k)!}$}x^{2k}}
\!-\!
{\Big (}\mbox{\small $\displaystyle \frac{2}{\pi}$}{\Big )}^{\!\!2m}\!\!
\left(\!\ln \frac{2}{\pi}
\!-\!
\displaystyle\sum\limits_{k=1}^{m-1}{
\mbox{\small $\displaystyle\frac{|B_{2k}| \pi^{2k}}{2 \,k\, (2k)!}$}
}\!\right)
\!x^{2m}.
\end{array}
$$

Let us denote by $H_{n,m}(x)$ the polynomial on the right-hand side of the above inequality.
Thus, for $x \!\in\! \left(0,\mbox{\small $\displaystyle\frac{\pi}{2}$}\right]$ and $m,n \!\in\! \NN, \,\, m,n \geq 2$ we have:
\begin{equation}
\label{H-polynomial-C2}
 F_{2}(x)  \, < \, H_{n,m}(x).
\end{equation}
Hence,  to prove  inequalitiy  (\ref{ekviv-C2-ineq})  it is sufficient to prove
 the following polynomial inequality:
\begin{equation}
\label{polynomial-C2}
H_{n,m}(x) \, < \, 0
\end{equation}
for $x \!\in\! \left(0,\mbox{\small $\displaystyle\frac{\pi}{2}$}\right]$, $n, m \!\in\! \NN, \, m\geq 2$.

\medskip
Let us consider the polynomial $ H_{n,m}(x)$ for $n=3$ and $m=3$.

\medskip
For the polynomial $H_{3,3}(x)$ it is not difficult to find  its smallest positive root  $x_1=1.0959152...$,
and  to determine the sign of   $H_{3,3}(x)$ for  $x \!\in\! \left(0, x_1\right)$:
$$
H_{3,3}(x)< 0.
$$

Thus,   inequality (\ref{polynomial-C2})  holds
for $n=3$, $m=3$ and every $x\in \left(0, x_1\right)$. Therefore,  inequality
(\ref{ekviv-C2-ineq})  holds true for every $x\in \left(0, x_1\right)$.\\

 Further, our goal is to prove  inequality (\ref{ekviv-C2-ineq})
 for every  $x \!\in\! \left[x_1,\mbox{\small $\displaystyle\frac{\pi}{2}$}\right]$.
 By   the following change of variables:
 $$\displaystyle x= \dfrac{\pi}{2}-t   ~~~~~ 
 $$
 for  $x \!\in\! \left[x_1,\mbox{\small $\displaystyle\frac{\pi}{2}$}\right]$ inequality  (\ref{ekviv-C2-ineq})
becomes
\begin{equation}\label{smena-C2-ineq}
F_2\left(\mbox{\small $\displaystyle\frac{\pi}{2}$}-t\right) < 0
\end{equation}
for $t \!\in\! \left[0, \dfrac{\pi}{2}-x_1\right]$.
We prove  inequality (\ref{smena-C2-ineq}) for
$$
t \!\in\!  \left[0,\, c \right)
$$
where  $ \dfrac{\pi}{2}\!-\!x_1 \!<\! c \!<\! \dfrac{\pi}{2}$.   Let us, for example, select $c=0.6$.\\

\smallskip
We have:
$$
\!\!\!\!\!
\begin{array}{rl}
 F_2\!\left( \mbox{\small $\displaystyle\frac{2}{\pi}$}\!-\!t\right)\! = \! & \!\!\!
\displaystyle P_{3}\!\left(\mbox{\small $\displaystyle\frac{2}{\pi}$}\!-\!t\right)
\ln\!\left(1\!-\!\mbox{\small $\displaystyle\frac{4(\pi-2)}{\pi^3}$}\left( \mbox{\small $\displaystyle\frac{2}{\pi}$}\!-\!t\right)^{\!2}\right)  -
\ln \cos {t}  -  \ln \left( \mbox{\small $\displaystyle\frac{2}{\pi}$}-t  \right)
\\ [2.1em]
 = & \!\!\!
\displaystyle P_{3}\!\left(\mbox{\small $\displaystyle\frac{2}{\pi}$}\!-\!t\right)
\!\left( \ln \! \mbox{\small $\displaystyle\frac{2}{\pi}$}
\!+\!
\ln \! \left(\! 1 \!+\! \mbox{\small $\displaystyle\frac{2(\pi-2)t(\pi-t)}{\pi^2}$}\!\right)\!\right)
 \!-\!  \ln \cos {t}
 \!-\!
 \ln\!\mbox{\small $\displaystyle\frac{\pi}{2}$}
\!+\!
\ln\left(1\!-\!\mbox{\small $\displaystyle\frac{2t}{\pi}$}\right).
 \end{array}
$$
From (\ref{Polinom-P3}) we conclude:
$$
P_{3}\!\left(\mbox{\small $\displaystyle\frac{\pi}{2}$}-t\right)
\! =\!
\alpha_2 \!\left(\mbox{\small $\displaystyle\frac{\pi}{2}$}-t\right)^{\!3}
+
\beta_2 \!\left(\mbox{\small $\displaystyle\frac{\pi}{2}$}-t\right)^{\!2}
+
\delta_2 >0
$$
for every $ t \!\in\! \left[0,  c \right)$. \\
As
$$\frac{2(\pi-2)t(\pi-t)}{\pi^2} < 1   ~~~~  \mbox{and} ~~~~ \frac{2t}{\pi} < 1$$
for every $t \!\in\! \left[0, \dfrac{\pi}{2}\right)$ i.e. for every  $t \!\in\! \left[ 0,  c \right)$,
and based on  (\ref{Procena-ln-cos-x}),  (\ref{Procena-ln(1+x)})   and (\ref{procena-ln(1-x)}),  we can conclude that
for $t \!\in\! \left[0, c\right)$ and  $ n\!=\!2l\!-\!1,  \, l \!\in\! \NN, \, m_1, m_2 \!\in\! \NN,  \, m_1 \geq 2 $
the following inequalities hold:
$$
\!\!\!\!\!
-\ln \cos t
\leq
\!\!\displaystyle\sum\limits_{k=1}^{m_1-1}{\!
\mbox{\small $\displaystyle\frac{2^{2k-1}(2^{2k}\!-\!1)|B_{2k}|}{k(2k)!}$}t^{2k}}
-
{\Big (}\mbox{\small $\displaystyle\frac{1}{c}$}{\Big )}^{\!\!2m_1}\!\!
\!\left(\!\ln{\cos{c}} - \!\!
\displaystyle\sum\limits_{k=1}^{m_1-1}{
\mbox{\small $\displaystyle\frac{2^{2k-1}(2^{2k}\!-\!1)|B_{2k}|}{k(2k)!}$}
c^{2k}}\!\!\right)\!t^{2m_1},
$$
$$\ln \! \left(\! 1 + \mbox{\small $\displaystyle\frac{2(\pi-2)t(\pi-t)}{\pi^2}$}\!\right)
\leq \displaystyle\sum\limits_{k=1}^{n}{\displaystyle (-1)^{k-1}\frac{2^k (\pi-2)^k}{k \pi^{2k}} t^k (\pi-t)^k
},
$$
$$
\ln\!\left(1-\dfrac{2t}{\pi}\right) \, \leq \,
 -\displaystyle\sum\limits_{k=1}^{m_2}{\displaystyle\frac{2^{k}t^{k}}{k\,\pi^k}}.
$$
Finally, based on the above results we have:
$$
\!\!\!\!\!\!\!\!\!\!\!\!
\begin{array}{lcl} F_2\!\left(\mbox{\small $\displaystyle\frac{\pi}{2}$}\!-\! t\right)
\!\!\!&\!\!<\!\!&\!\!\!
\,\displaystyle P_{3}\!\left(\mbox{\small $\displaystyle\frac{\pi}{2}$}\!-\!t\right)
\!\!\left( \ln \! \mbox{\small $\displaystyle\frac{\pi}{2}$}
+
\displaystyle\sum\limits_{k=1}^{n}{\displaystyle (-1)^{k-1}\frac{2^k (\pi-2)^k}{k\, \pi^{2k}} t^k \,(\pi-t)^k
}
\right) \, +  \\  [1.5em]
\!\!\!&\!\!+\!\!&\!\!\!
\displaystyle\sum\limits_{k=1}^{m_1-1}{\!
\mbox{\small $\displaystyle\frac{2^{2k-1}(2^{2k}\!-\!1)|B_{2k}|}{k(2k)!}$}t^{2k}}
\!-\!
{\Big (}\mbox{\small $\displaystyle\frac{1}{c}$}{\Big )}^{\!\!2m_1}\!\!
\left(\ln{\cos{c}} - \!\!\!
\displaystyle\sum\limits_{k=1}^{m_1-1}{\!\mbox{\small $\displaystyle\frac{2^{2k-1}(2^{2k}\!-\!1)|B_{2k}|}{k(2k)!}$}
 c^{2k}}\!\right)
\!t^{2m_1}  \\ [1.5em]
\!\!&\!\!-\!\!&\!\!
 \ln\!\mbox{\small $\displaystyle\frac{\pi}{2}$}
 -  \displaystyle\sum\limits_{k=1}^{m_2}{\displaystyle\frac{2^{k}t^{k}}{k\,\pi^k}}\,.
\end{array}
$$
Let us denote by $T_{n,m_1,m_2}(t)$ the polynomial on the right-hand side of the above inequality.
Thus, for $t \!\in\! \left[0, c\right)$ we have:
\begin{equation}
\label{T-polynomial-C2}
 F_{2}\!\left(\dfrac{\pi}{2}-t\right)  \, < \, T_{n,m_1,m_2}(t).
\end{equation}

Hence, for the  proof of  inequality  (\ref{smena-C2-ineq}) it is sufficient to prove
 the following polynomial inequality:

\begin{equation}
\label{polynomial-C2-drugi}
T_{n,m_1,m_2}(t) \, < \, 0
\end{equation}
for every $t \!\in\! \left[0,\, c\right)$,   $n,m_1, m_2 \!\in\! N, \, m_1 \geq 2 $.

\smallskip
Let us consider the polynomial $T_{n,m_1, m_2}(t) $ for $n=5$, $m_1=7$  and $m_2=9$.

\smallskip
For the polynomial $T_{5,7,9}(t)$ it is not difficult to find  its smallest positive
root  $t_1=0.6257524\ldots$,
and to determine the sign of  $T_{5,7,9}(t)$ for  $t \!\in\! \left[0, t_1\right)$:
$$T_{5,7,9}(t)< 0. $$

Thus,  inequality (\ref{polynomial-C2-drugi})  holds
for $n=5, m_1=7,  m_2=9$ and every \mbox{$t\in \left[0, c\right)$.}
Therefore,  inequality
 (\ref{smena-C2-ineq}) holds true for every
 $t \!\in\! \left[0,\mbox{\small $\displaystyle\frac{\pi}{2}$}-x_1\right] \!\subset\! \left[0,\, c\right)$,  and
 ine\-quality (\ref{ekviv-C2-ineq}) holds for  every  $x \!\in\! \left[x_1,\mbox{\small $\displaystyle\frac{\pi}{2}$}\right)$.

\smallskip
The proof of Theorem 2 is now complete. \hfill $\Box$

\bigskip
The following theorem,  stated as a conjecture in \cite{Nishizawa_2015}, was proved in \cite{Malesevic_Lutovac_Banjac_2018} using the
approximations and methods  from \cite{Mortici_2011}, \cite{Malesevic_Makragic_2016}, \cite{Nenezic_Malesevic_Mortici_2016},
\cite{Lutovac_Malesevic_Mortici_2017}, \cite{Milica_Makragic_2017} and \cite{Malesevic_Lutovac_Rasajski_Mortici_2017}.
In this paper we give another proof of this conjecture. In particular, we show that this conjecture is a consequence of Theorem 2.

\begin{theorem}
Let the function
$$
\displaystyle
f_3(x)
=
\left(1 - \mbox{\small $\displaystyle\frac{4(\pi-2)}{\pi^3}$}\,x^2\right)^{\alpha_3 x^3  + \delta_3}
\!-\!
\frac{\sin x}{x}
$$
for $x \!\in\! \left(0,\mbox{\small $\displaystyle\frac{\pi}{2}$}\right]$ satisfy  the following conditions:
\begin{equation}
\label{Uslovi_3}
f_{3}(0+)
\!=\!
f_{3}'(0+)
\!=\!
f_{3}''(0+)
\!=\!
0,\;
f_{3}\!\left(\mbox{\small $\displaystyle\frac{\pi}{2}$}\right)
\!=\!
0.
\end{equation}
Then:
\begin{equation}
\label{Koeficijenti_3}
\begin{array}{l}
\alpha_3 = -\mbox{\small $\displaystyle\frac{\pi^3-24\pi+48}{3(\pi-2)\pi^3}$}, \\[1.5 ex]
\delta_3 =  \mbox{\small $\displaystyle\frac{\pi^3}{24(\pi-2)}$}.
\end{array}
\end{equation}
and
\begin{equation}
f_3(x)<0,
\end{equation}
for every $x \!\in\! \left(0,\mbox{\small $\displaystyle\frac{\pi}{2}$}\right].$
\end{theorem}

~\\
\noindent {\bf Proof.}
As shown in  Subsection~\ref{About-constraints},
the conditions in (\ref{Uslovi_3}) yield a  system   of linear equations  (shown in ~(\ref{fj''(0)=0}) and
(\ref{fj(Pi/2)=0}))  in variables $\, \alpha_3$ and $ \delta_3$.
The symbolic values (\ref{Koeficijenti_3}) of  $\alpha_3$ and $ \delta_3$ are obtained by solving this system.

 Based on Theorem 2, it is enough to prove that for every $x \!\in\! \left(0, \dfrac{\pi}{2}\right]$
 $$
 \left(1 - \mbox{\small $\displaystyle\frac{4(\pi-2)}{\pi^3}$}\,x^2\right)^{\alpha_3 x^3  + \delta_3}
\!<
\left(\!1 - \mbox{\small $\displaystyle\frac{4(\pi-2)}{\pi^3}$}\,x^2\!\right)^{\alpha_2 x^3 + \beta_2 x^2  +
\delta_2}.
$$
The above inequality is equivalent to the following inequalities:
$$
\alpha_3 x^3  + \delta_3  >  \alpha_2 x^3 + \beta_2 x^2  +  \delta_2
\;\;\, \Longleftrightarrow \;\;\,
x^2\left((\alpha_3  -  \alpha_2) x - \beta_2 \right)  > 0.
$$
It is not hard to check that the  polynomial inequality on the right-hand side holds true for every
$x \!\in\! \left(0, \mbox{\small $\displaystyle\frac{\pi}{2}$}\right]$.
The proof of Theorem 3 is therefore complete.

\hfill $\Box$

\section{Conclusion}

In this paper we presented a new approach to proving some exponential inequalities.
We illustrated our ideas in the proofs of Theorems 1 and 2, i.e. in the proofs
of some exponential inequalities connected with the sinc function. Using particular
approximations based on the power series expansions, proving of the exponential inequalities
was reduced to proving of the corresponding polynomial inequalities.

\smallskip
Our approach can be applied more broadly than just to exponential inequalities.
Another potential application of our technique is in the area of establishing
new polynomial bounds as shown, for example, in inequalities (\ref{G-polynomial-C1})
and (\ref{H-polynomial-C2}) in the proofs of Theorem 1 and 2, respectively.

\bigskip
\noindent
\textbf{Competing Interests.} The authors would like to state that they do not have any competing interests in the subject
of this research.

\bigskip
\noindent
\textbf{Author's Contributions.} All the authors participated in every phase of the research conducted for this paper.

\break

\smallskip


\end{document}